\documentclass[11pt]{amsart}

\usepackage[colorlinks=true, pdfstartview=FitV, linkcolor=blue, 
      citecolor=blue]{hyperref}
\usepackage{amsfonts, amssymb, amsmath, mdframed, graphicx}
\usepackage{bbm, a4wide}

\DeclareMathAlphabet{\mymathbb}{U}{bbold}{m}{n}

\newcommand{\tdens}{\tau}
\newcommand{\edens}{\varepsilon}

\newcommand{\be}{\begin{equation}}
\newcommand{\ee}{\end{equation}}

\newtheorem{theorem}{Theorem}

\theoremstyle{definition}

\newtheorem{remark}[theorem]{Remark}

\newcommand{\ER}{Erd\H{o}s-R\'enyi}

\title{Tripodal structure in undersaturated random graphs}

\author{William DiCarlo}
\address{Department of Mathematics, University of Texas, \newline
  \indent 2515 Speedway, PMA 8.100 Austin, TX 78712, USA}
\email{will.a.dicarlo@utmail.utexas.edu}

\author{Lorenzo Sadun}
\address{Department of Mathematics, University of Texas, \newline
  \indent 2515 Speedway, PMA 8.100 Austin, TX 78712, USA}
\email{sadun@math.utexas.edu}

\begin{document}

\begin{abstract}
We numerically investigate typical graphs in a region of the Strauss model of random 
graphs with constraints on the densities of edges and triangles. This region, where typical graphs had been expected
to be bipodal but turned out to be tripodal, involves edge densities $e$ below $e_0 = (3-\sqrt{3})/6 \approx 0.2113$ 
and triangle densities $t$ slightly below $e^3$. We determine the extent of this region in $(e,t)$ space 
and show that there is a discontinuous phase transition at the boundary between this region and a bipodal phase. 
We further show that there is at least one phase transition within this region, where the parameters 
describing typical graphs change discontinuously. 
\end{abstract}


\maketitle


\section{Introduction and results}\label{sec:intro}

\subsection{Background}

We are investigating the structure of large random graphs with specified densities $(e,t)$ of edges and 
triangles, the so-called Strauss model \cite{Strauss}. More precisely, we consider the ensemble of all graphs on $n$
vertices that have the specified densities of edges and triangles, up to a small tolerance that goes to 
zero as $n \to \infty$. For most values 
of $(e,t)$, there is a typical structure to such graphs, with all but a tiny fraction  
having essentially the same statistical properties. 
The fundamental question is understanding how this structure varies with $(e,t)$. In this paper we study this
question when $t$ is slightly below $e^3$ and $e$ is less than the threshold $e_0=(3-\sqrt{3})/6 \approx 0.2113$.

The key algebraic tool is {\em graphons}. (See \cite{Lov} for an overview of this expansive subject and 
\cite{BCL, BCLSV, LS1} for some original references.) 
A graphon is a function $g: [0,1]^2 \to [0,1]$ with $g(x,y)=g(y,x)$. This can be viewed
either as a limit of large graphs or as a random process for producing graphs. To generate a graph on $n$
vertices from a graphon, first pick $n$ numbers $x_i$ independently and uniformly on $[0,1]$. Then flip $\binom{n}{2}$ independent
biased coins, assigning an edge to vertices $i$ and $j$ with probability $g(x_i, x_j)$. The statistical properties
of graphs that are generated this way are easy to read off from the graphon. For instance, the expected
densities of edges and triangles are
\begin{eqnarray}
\edens(g) &:=& \int_0^1 \! \int_0^1 g(x,y) \, dx \, dy \qquad \qquad  \hbox{and} \cr 
\tdens(g) &:=& \int_0^1 \! \int_0^1 \! \int_0^1 g(x,y) g(y,z) g(z,x) \, dx \, dy \, dz, 
\end{eqnarray}
respectively. As $n \to \infty$, the variance in the edge and triangle densities go to zero, so almost all graphs generated
by $g$ have edge and triangle densities close to $(\edens(g), \tdens(g))$. We can thus speak of $\edens(g)$ and $\tdens(g)$
as the edge/triangle densities of the graphon, and not just as the {\em expected} densities. 

The number of graphs that are created
in this way is determined by the Shannon entropy of the graphon, 
\begin{equation} \label{eq:S(g)} S(g) := \int_0^1 \! \int_0^1 H(g(x,y)) \, dx \, dy, \end{equation}
where
\[ H(u) = - (u \ln(u) + (1-u) \ln(1-u)) \]
is the usual entropy of a coin flip with probability $u$ of getting heads. 
Thanks to the large deviations principle of Chatterjee and Varadhan \cite{CV} and to subsequent work by Radin
and Sadun \cite{RS1}, understanding the typical structure of graphs with edge/triangle densities $(e,t)$ is 
equivalent to finding the graphon $g$ that maximizes $S(g)$ subject to the constraints 
$(\edens(g), \tdens(g))=(e,t)$. 
If $\tilde g(x,y) = g(\sigma(x), \sigma(y))$, where $\sigma$ is a measure-preserving transformation
of $[0,1]$, then $\tilde g$ and $g$ have the same densities and the same entropy. Indeed, they generate
the same graphs on $n$ vertices with exactly the same probabilities! When we apply constraints, maximize the entropy, 
and speak of the optimal graphon being unique, we always mean ``unique up to measure-preserving transformation.'' 

If the entropy-maximizing graphon $g$ with edge/triangle densities $(e,t)$ 
is unique (up to measure-preserving transformations), then the ensemble of large 
graphs with edge/triangle densities close to $(e,t)$ is essentially the same as the ensemble of large graphs generated
by $g$. Specifically, for any $\epsilon>0$ there exists a constant $K>0$ such that all but a fraction $e^{-Kn^2}$
of the graphs on $n$ vertices with edge/triangle densities close to $(e,t)$ are within $\epsilon$ of $g$ in the 
``cut metric'' \cite{RS1, RS3}. 
Conversely, all but an exponentially small fraction of the graphs generated by $g$ have edge/triangle densities close to 
$(e,t)$. The combinatorial problem of understanding typical graphs with specified edge/triangle densities thus reduces to 
the analytic problem of maximizing $S(g)$ subject to the constraints $(\edens(g), \tdens(g))=(e,t)$.

\begin{figure}[ht]
\includegraphics[width=4in]{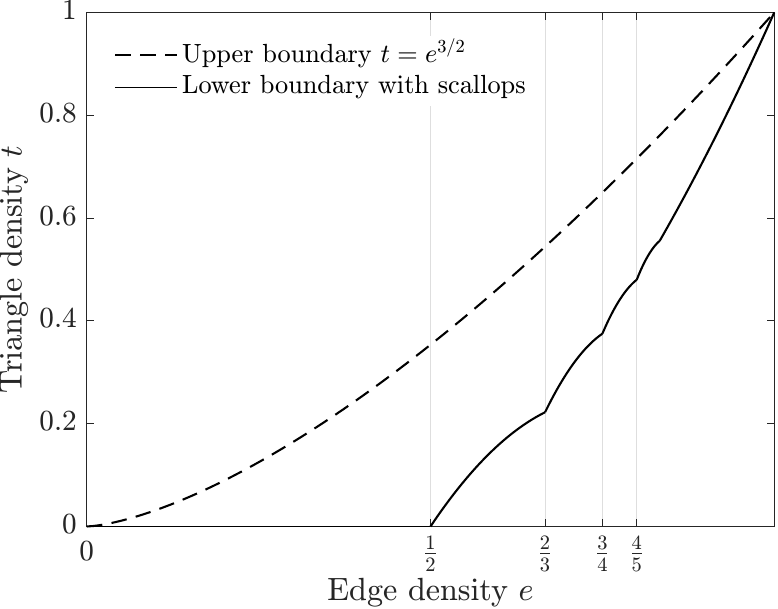}
\caption{The Razborov triangle. The curvature of the ``scallops'' on the lower 
right is exaggerated for visibility.}
\label{fig:Razborov}
\end{figure}

The space of all possible values of $(e,t)$ in the edge-triangle model 
was first worked out by Razborov \cite{Raz} (see also \cite{PR}) and is called the Razborov triangle. See Figure \ref{fig:Razborov}. 
A {\em phase} is an open subset of the Razborov triangle where the typical structure is a real-analytic function
of $(e,t)$ in the sense that
the optimal graphon is unique and depends analytically on $(e,t)$. 
This implies that the density of any fixed subgraph, such as a square or pentagon or tetrahedron, is also an 
analytic function of $(e,t)$. 
At a phase boundary, the optimal graphon fails to be analytic and may even change discontinuously, implying
an abrupt change in the density of some subgraph. 

\begin{figure}[ht]
\includegraphics[width=4in]{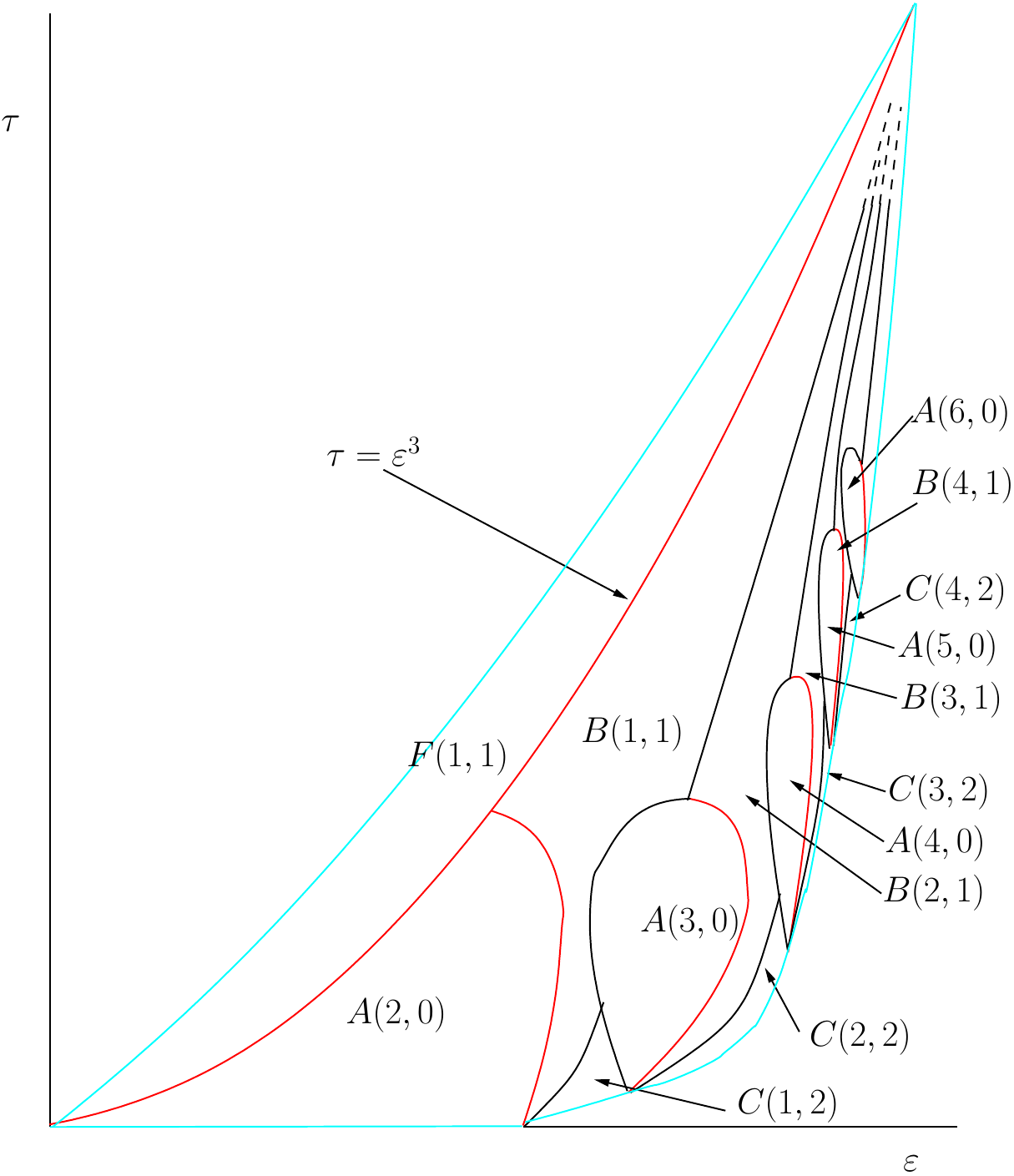}
\caption{A conjectured phase portrait in the edge/triangle model. The size and shape of 
each phase are only schematic. In a realistic rendering, many of the phases would be too small to see.}
\label{fig:conjecture}
\end{figure}

In 2017, \cite{KRRS1} conjectured that the Razborov triangle consisted of one phase above the ``\ER'' curve $t=e^3$ and
three infinite families of phases below the \ER{} curve. See Figure \ref{fig:conjecture} for a schematic phase diagram.
The conjectured phase picture has been proven to be correct
\begin{itemize} 
\item Just above the \ER{} curve. \cite{KRRS3}
\item Just below the \ER{} curve when $e>1/2$. \cite{NRS1}
\item Near the line segment $e=1/2$, $0<t<1/8$. \cite{NRS3}
\item Just below the top boundary $t=e^{3/2}$. \cite{RS2}
\item Just above the bottom boundary when $e<1/2$. \cite{RS2}
\item Just above each of the ``scallops'' on the lower right portion of the boundary, with the $k$-th scallop 
running from $e=\frac{k}{k+1}$ to $e = \frac{k+1}{k+2}$. \cite{RS2} 
\end{itemize} 
That is, the optimal graphons in 
open subsets of the $F(1,1)$, $B(1,1)$, $A(2,0)$, and $C(k,2)$ phases, for
every positive integer $k$, have been proven to exist, to be unique, and to have the form conjectured in \cite{KRRS1}.

However, there is also a region, with $e < e_0 := (3-\sqrt{3})/6 \approx 0.2113$ and with $t$ slightly below
$e^3$, where the optimal graphon is {\bf not} described by the $A(2,0)$ phase of Figure \ref{fig:conjecture} \cite{NRS3}! Unfortunately, 
the proof that the typical structure is not $A(2,0)$ did not provide any estimates for how big 
this exceptional
region is or what the typical structure is. \cite{NRS3} exhibited a ``tripodal'' graphon that has more entropy 
than the ``bipodal'' $A(2,0)$ graphon, but did not optimize the parameters of this tripodal graphon or 
rule out the possibility that there is a different kind of graphon with even more entropy.

In other words, Figure \ref{fig:conjecture} needs to be amended, adding one or more mystery phases just below 
the \ER{} curve for $e < e_0$. 
In this paper we report on some numerical explorations of this {\em terra incognita}.

\subsection{Multipodality and the ansatz of \cite{NRS3}}

A graphon $g$ is said to be $k$-podal if we can break the unit interval into $k$ disjoint measurable regions
$I_1, \ldots, I_k$ such that $g(x,y)$ is constant on each $I_i \times I_j$. A graphon that is 
$k$-podal for some $k$ is said to be {\em multipodal}. The words ``bipodal'' and ``tripodal'' mean 2-podal and
3-podal, respectively. The regions $I_i$ are called ``podes''. After applying a measure-preserving transformation
of $[0,1]$, we can assume that each pode $I_i$ is an interval and each $I_i \times I_j$ is a rectangle. 
A $k$-podal graphon $g$ is then described by the sizes $|I_i|$
of the podes and by a symmetric $k\times k$ matrix giving the values of $g$ on each rectangle. 

An $(n+m)$-podal graphon is said to have {\em $(n,m)$ symmetry} if it is invariant under permutation of the first $n$
podes and also invariant under permutation of the last $m$ podes. In the conjectured phase diagram of Figure 
\ref{fig:conjecture}, all phases are multipodal with either $(n,0)$, $(n,1)$ or $(n,2)$ 
symmetry, as indicated by the name of each phase. (In related models involving $k$-stars, all entropy maximizers have been proven to be multipodal \cite{KRRS2}.
However, proving multipodality for all phases in the edge/triangle model remains an open problem.) 

\begin{figure}[ht]
\includegraphics[width=3in]{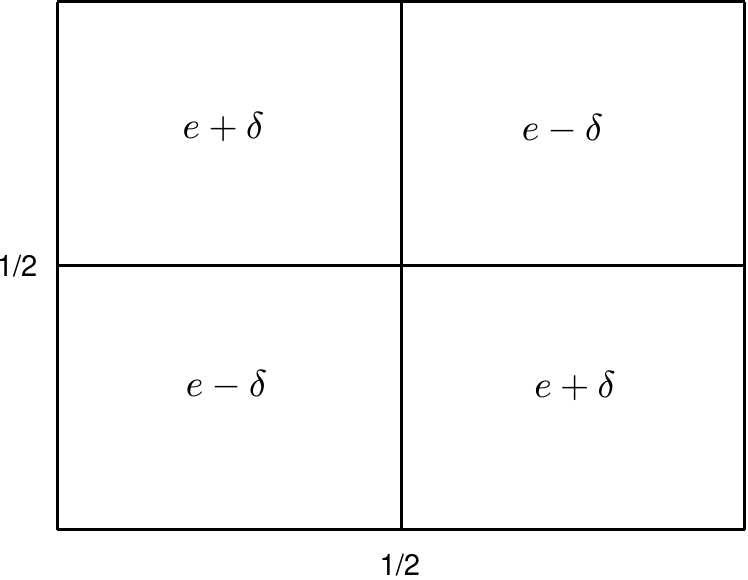}
\caption{The unique symmetric bipodal graphon with edge and triangle densities $(e, e^3-\delta^3)$.}
\label{fig:Bipodal}
\end{figure}

The space of $(2,0)$-symmetric graphons is only 2-dimensional, so there is a unique graphon with edge/triangle
densities $(e,t)$. This is shown in Figure \ref{fig:Bipodal}, where $t=e^3-\delta^3$. 
The entropy of this ``symmetric bipodal'' graphon $g_{sb}$ is 
\be \label{eq:Ssb} 
S(g_{sb}) = \frac12 \Big(H(e+\delta) + H(E-\delta)\Big) = H(e) + \frac12 H''(e) \delta^2 + O(\delta^4). \ee

\begin{figure}[ht]
\includegraphics[width=3in]{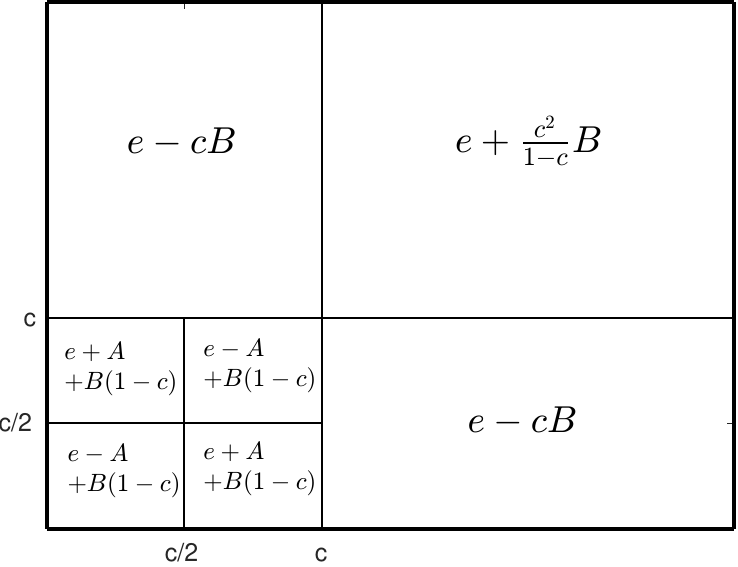}
\caption{A tripodal graphon with (2,1) symmetry and constant degrees.}
\label{fig:ansatz}
\end{figure}

Now consider the $(2,1)$-symmetric graphon of Figure \ref{fig:ansatz}, 
depending on parameters $A$, $B$, and $c$. 
Specifically, we have podes $I_1$ and $I_2$, each of size $c/2$, and $I_3$ of size $1-c$, and our graphon
takes on the values
\be \label{eq:ansatz} g(x,y) = \begin{cases} e-A+B(1-c) & (x,y) \in (I_1 \times I_1) \cup (I_2 \times I_2), \cr 
e+A +B(1-c) & (x,y) \in (I_1 \times I_2) \cup (I_2 \times I_1), \cr 
e-cB & (x,y) \in (I_1 \times I_3) \cup (I_2 \times I_3) \cup (I_3 \times I_1) \cup (I_3 \times I_2), \cr 
e + \frac{c^2}{1-c} B & (x,y) \in I_3\times I_3. \end{cases} \ee
The edge density is $e$ and the triangle density is $e^3 - c^3 (A^3-B^3)$, so $c=\delta (A^3-B^3)^{-1/3}$. 
Holding $A$ and $B$ fixed and considering the behavior as $c \to 0$, we have 
\begin{eqnarray} S(g) &=& \frac{c^2}{2} \Big ( H(e-A+B(1-c)) + H(e+A+B(1-c))\Big ) \cr 
&& + 2c(1-c) H(e-cB) + (1-c)^2 H\left ( e + \frac{c^2}{1-c} B \right ) \cr 
& = & H(e) + \frac12 F(A,B) \delta^2 + O(\delta^3),
\end{eqnarray}
where
\be \label{eq:F} F(A,B) = \frac{H(e+A+B) + H(e-A+B) -2H(e) -2BH'(e)}{(A^3-B^3)^{2/3}} \ee
and we have used Taylor series to approximate $H(u)$ when $u \approx e$. 
If we can find parameters $A$ and $B$ such that $F(A,B)> H''(e)$, then for sufficiently small $\delta$ 
we have $S(g) > S(g_{sb})$. In \cite{NRS3}, a power series expansion was used to show
that such values of $(A,B)$ exist, with $A$ small and $B=O(A^2)$, whenever $e < e_0 := (3-\sqrt{3})/6$. 

\begin{remark}
The expression (\ref{eq:ansatz}) may look complicated and unintuitive. However, it is actually one of the simplest 
possible structures to improve on (\ref{eq:Ssb}). By integrating a simple Taylor series expansion of $H(g(x,y))$ 
around $g(x,y)=e$ we obtain
\be \label{eq:S-expand} 
S(g) = H(e) + \frac12 H''(e) \|g-e\|_{L^2}^2 + O(\|g-e\|_{L^2}^3) \le H(e) + \frac12 H''(e) \delta^2 
+ o(\delta^2), \ee
for all graphons $g$ that are pointwise close to $e$. In order to obtain a coefficient of $\delta^2$ that
is greater than $\frac12 H''(e)$, we must consider instead a graphon that differs greatly from $e$ on a subset 
of the unit square of measure $O(\delta^2)$. 

\begin{figure}[ht]
\includegraphics[width=3in]{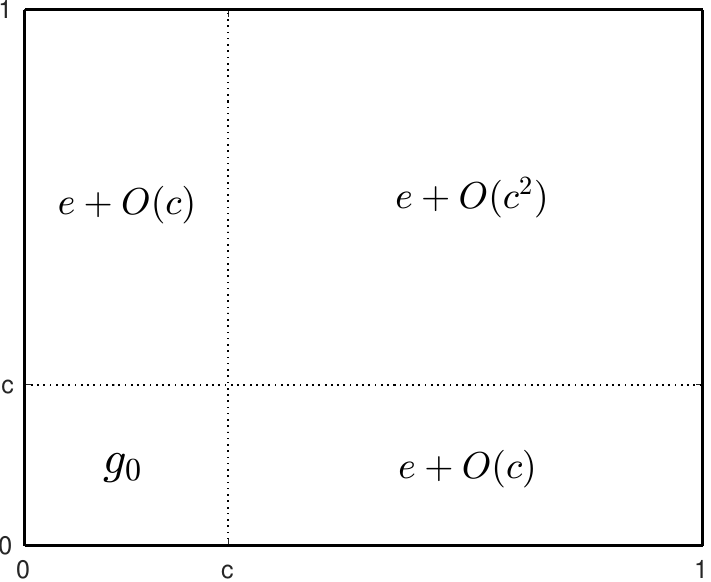}
\caption{A strategy for generating efficient graphons. The $O(c)$ and $O(c^2)$ terms are chosen to make the 
degree function exactly constant.}
\label{fig:graphon-g0}
\end{figure} 

This suggests the following general construction. See Figure \ref{fig:graphon-g0}. Start with an arbitrary graphon $g_0$, with no constraints on the
edge or triangle densities of $g_0$. Let $I_1=[0,c)$ for some small constant $c$ and let $I_2 = [c,1]$. Let 
\be \label{eq:ansatz2} g(x,y) = \begin{cases} g_0\left (\frac{x}{c}, \frac{y}{c}\right ) & (x,y) \in I_1 \times I_1, \cr 
e - \frac{c}{1-c} d(x/c) & (x,y) \in I_1 \times I_2, \cr 
e - \frac{c}{1-c} d(y/c) & (x,y) \in I_2 \times I_1, \cr 
e + \frac{c^2}{(1-c)^2} B & (x,y) \in I_2 \times I_1. \end{cases} \ee
where $d(x)=\int_0^1 g_0(x,y) \, dy$ is the degree function of $g_0$ 
and $B = \int_0^1 (d(x)-e) \, dx = \iint (g_0(x,y)-e) \, dx \, dy $ is the 
difference between the edge density of $g_0$ and $e$. That is, we put a rescaled copy of $g_0$ in the lower left
corner of the unit square and adjust $g$ on the rest of the square to keep the degree function constant, 
thereby minimizing
$\|g-e\|_{L^2}$ for given $\delta$ and $g_0$. 

The triangle density for small $c$ is easily computed in terms of $g_0$: 
\be
\tdens(g) =  e^3 + c^3 \tdens(g_0-e) + O(c^4), \ee
so 
\be 
\delta =  c \, \tdens(e-g_0)^{1/3} + O(c^2). 
\ee 
We similarly compute the entropy of $g$ to leading order in $c$, obtaining
\begin{eqnarray} \label{eq:F-general}
S(g) & = & c^2 S(g_0) + (1-c^2) H(e) - c^2 B H'(e) + O(c^3) \cr 
& = & H(e) + c^2 (S(g_0) - H(e) - B H'(e)) + O(c^3) \cr 
& = & H(e) + \frac{\delta^2}{2} \left ( \frac{2 S(g_0) - 2H(e) - 2B H'(e)}{\tdens(e-g_0)^{2/3}} \right ) + O(\delta^3). 
\end{eqnarray} 

The simplest version of this construction involves a constant graphon $g_0$. That is {\em precisely} the asymptotic
form of the optimal graphon when $e>1/2$ and $t$ is slightly less than $e^3$ \cite{NRS1}. 

The next simplest possibility is for $g_0$ to be symmetric bipodal. In that case,
the general construction (\ref{eq:ansatz2}) reduces to (\ref{eq:ansatz}), albeit
with the constant $B$ rescaled by a factor of $1-c$. 
The coefficient
of $\frac{\delta^2}{2}$ in (\ref{eq:F-general}) reduces to the function $F(A,B)$ of (\ref{eq:F}).  

\end{remark}

\subsection{Organization of the paper}

We extend the analysis of \cite{NRS3} with two goals in mind. First, we wish to find the actual entropy 
maximizing graphon $g$, not just a graphon that does better than $g_{sb}$. Second, we wish to explore the
entire tripodal region, not just an infinitesimal neighborhood of the \ER{} curve.

\begin{enumerate}
\item In Section \ref{sec:F} we determine the values of $A$ and $B$ that maximize $F(A,B)$, over a range of values of $e$. This is a problem
in 2-variable calculus, not in functional analysis, and does not require any fancy methods. A combination of
sampling on a grid and Newton's method is sufficient. As expected, it is possible to do
better than $H''(e)$ precisely when $e < e_0$. We also discover a transition at $e \approx 0.0024$, where the best values 
of $A$ and $B$ change discontinuously. 

\item The set of possible graphons with (2,1) symmetry is slightly more general than the ansatz (\ref{eq:ansatz}),
with the degree function on the two podes differing in general by a quantity $D/2$, with the original 
ansatz corresponding to $D=0$. In Section \ref{sec:D} we explore (2,1)-symmetric graphons with $D$ nonzero and show 
that the optimal value of 
$D$ is $O(\delta^2)$, yielding an entropy that is $O(\delta^4)$ better than with $D=0$. Since we are mostly looking at the 
entropy at order $\delta^2$, this $O(\delta^4)$ contribution does not change the overall picture. 
For a qualitative understanding of the tripodal phase and of phase transitions,
the ansatz (\ref{eq:ansatz}) is enough. 

\item In Section \ref{sec:limbo}, we probe the size of the tripodal phase. Specifically, for each value of $e$ we
determine the range of $\delta$ values for which a (2,1)-symmetric tripodal graphon $g$ has more entropy than a 
symmetric bipodal graphon. 
When $\delta$ reaches a certain size, the entropies of $S(g)$ and $S(g_{sb})$ become equal, 
indicating a 
phase transition. The graphon $g$ does not become singular as we approach this point, nor does it approach
$g_{sb}$. Instead, the optimal graphon changes discontinuously from $g$ to $g_{sb}$. 

We do this analysis twice, once for the ansatz and once allowing $D$ to vary. There are small quantitative
differences, but the entropies are close and
the optimal values of $(A,B)$ are nearly the same. 

\item Finally, in Section \ref{sec:small-e}, we explore the transition near $e=0.0024$. 
Assuming that the optimal graphon is either symmetric bipodal or (2,1)-symmetric tripodal,
we determine the boundary between two distinct tripodal phases. 
\end{enumerate}

\section{The maximum of $F(A,B)$} \label{sec:F} 
The ansatz (\ref{eq:ansatz}) gives the most general graphon with (2,1) symmetry, with edge density $e$, and 
with constant degree function $d(x) = \int_0^1 g(x,y) \, dy$. 
Maximizing $F(A,B)$ is tantamount to maximizing the entropy among graphons with (2,1) symmetry and with triangle density just below $e^3$. 

\begin{figure}[ht]
\includegraphics[width=5in]{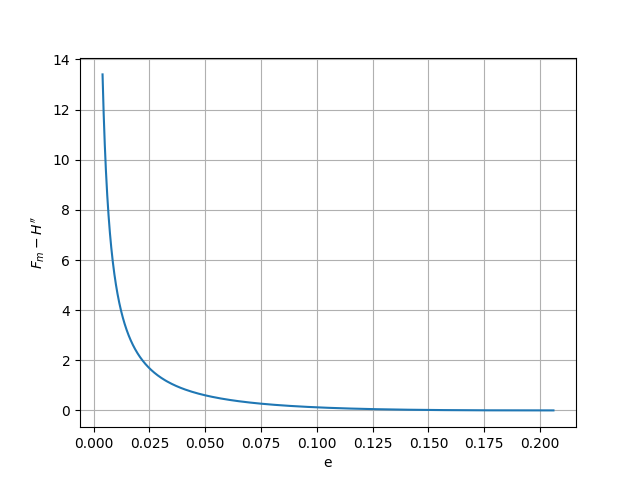}
\caption{$F_m(e)-H''(e)$ as a function of $e$.}\label{fig:F-ddHvse}
\end{figure} 

\begin{figure}[ht] 
\includegraphics[width=5in]{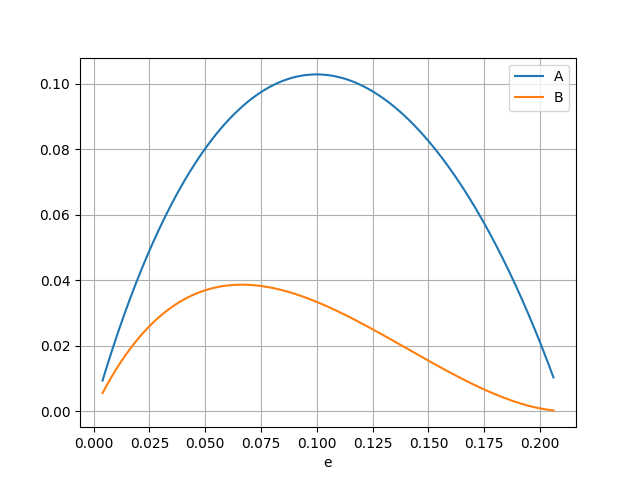}
\caption{The optimal values of $A$ and $B$ as functions of $e$.}\label{fig:ABvse}
\end{figure}

Let $F_m(e)$ be the maximum of $F(A,B)$ for a given edge density $e$. In order to compute $F_m(e)$, 
we first choose a value of $e$ and plot $F(A,B)$ in Desmos two ways, 
once for each variable while leaving the 
other variable on a slider. By moving the sliders and looking at the corresponding cross sections of $F(A,B)$, we
determine the approximate locations of the local maxima. We then apply Newton's method, starting at these
approximate locations, to localize the maximum to machine
accuracy. Having done this for one value of $e$, we increment $e$ in steps of $0.001$ and use the optimal values of 
$(A,B)$ from the previous $e$ as the starting point for Newton's method for the new value of $e$. This process 
successfully evaluated the optimal values of $A$ and $B$ and the values of $F_m(e)$ for $e$ between 
0.033 and 0.206. See Figures  \ref{fig:F-ddHvse} and 
\ref{fig:ABvse}. The same procedure, only with a smaller increment for $e$, allowed us to probe 
values of $e$ that are smaller than 0.033 or
larger than 0.206. 

\begin{figure}[ht]
\includegraphics[width=5in]{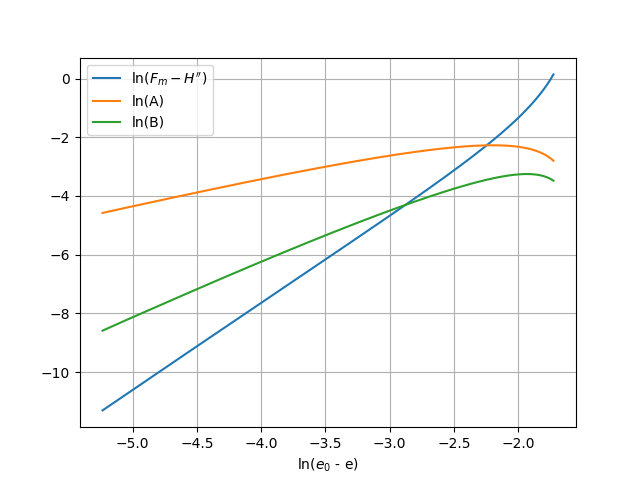}
\caption{The slopes of this log plot show how various functions scale with $e$ as $e \to e_0$. 
$A$ scales as $(e_0-e)^1$, $B$ as $(e_0-e)^2$, and $F_m-H''$ as $(e_0-e)^3$.} \label{fig:scaling}
\end{figure}

The relevant values of $A$ and $B$ go to zero as $e \to e_0$. To determine the rate of decay, we plot $\ln(A)$, $\ln(B)$ and $\ln(F_m-H'')$ versus $\ln(e_0-e)$ in Figure \ref{fig:scaling}. 
The plots are approximately linear as $e \to e_0$, with slopes 1, 2 and 3, respectively, indicating that 
$A \sim (e_0-e)^1$, $B \sim (e_0-e)^2$ and $F_m-H'' \sim (e_0-e)^3$.

Our procedure for maximizing $F(A,B)$ is guaranteed to find a {\em local\/} maximum near the previous best values of 
$(A,B)$. To ensure that we are 
not missing a separate (and potentially higher) local maximum, we re-examined the function $F(A,B)$ on Desmos by eye
for selected values of $e$. This turned out to be unnecessary for $e > 0.01$, where the local maximum found by our 
iterated Newton's method procedure is indeed the global maximum. However, things are different when $e$ is very small. 

\begin{figure}[ht] 
\includegraphics[width=5in]{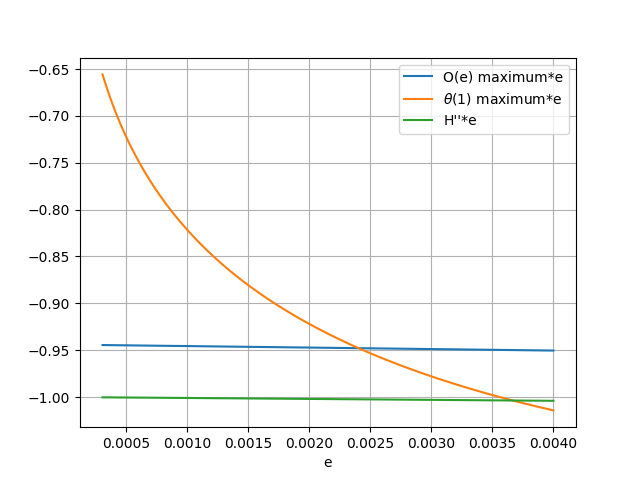}
\caption{The two local maxima of $F(A,B)$ for small $e$, compared to $H''(e)$.}\label{fig:Low-e}
\end{figure}

This is seen in Figure \ref{fig:Low-e}. There are two local maxima in the $F(A,B)$ plot. One, with $A \approx 2.5e$ and
$B \approx 1.5e$, yields values of $F(A,B)$ that are a few percent better than $H''(e)$, being close to $-0.95/e$ 
while $H''(e) \approx -1/e$. The other, with $A$ and $B$ close to 1/2, dominates when $e < 0.0024$. 
To get an algebraic understanding of this numerical result, consider $F(A,B)$ when $A=\frac12$ and 
$B=\frac12 - e$. (These particular values, while 
not optimal, give a lower bound for $F_m(e)$ and an upper bound for $F_m(e)/H''(e)$.)  In that case, 
\begin{eqnarray}
F(A,B) & = & \frac{H(e+A+B) + H(e-A+B) - 2H(e) - 2BH'(e)}{(A^3-B^3)^{2/3}} \cr 
& = & \frac{-2H(e) - (1-2e)H'(e)}{\left ( \frac34 e - \frac32 e^2 + e^3 \right )^{2/3}},
\end{eqnarray}
since $H(e+A+B)=H(1)=0$ and $H(e-A+B)=H(0)=0$. As $e \to 0$, the numerator scales as $\ln(e)$, while the 
denominator scales as $e^{2/3}$, so $F(A,B)$ scales as $e^{-2/3} \ln(e)$.
Since $H''(e) = -\left (\frac{1}{e} + \frac{1}{1-e} \right )$ scales as $e^{-1}$, the ratio
$F_m(e)/H''(e)$ goes to zero at least as fast as $e^{1/3} \ln(e)$.

The switch from one local maximum to another at $e \approx 0.0024$ indicates that there must be at least one phase transition
when $e$ is small, in addition to the transition between tripodal and bipodal graphons. However, it does not tell us what
the phase transition looks like. 

The simplest possibility is that there are two (2,1)-symmetric phases, with the boundary between them hitting the
\ER{} curve at $e \approx 0.0024$. Recall that maximizing $F(A,B)$ is tantamount to determining the best (2,1)-symmetric
graphon for $t$ slightly less than $e^3$. When $e \gtrsim 0.0024$ and $t$ is slightly less than $e^3$, the optimal graphon
is (2,1)-symmetric tripodal with all values on the order of $e$, in particular with all values less than 0.01. When 
$e \lesssim 0.0024$ and $t$ is slightly less than $e^3$, the optimal graphon takes values close
to 1 on $I_1 \times I_2$. On the curve separating the two phases, the optimal graphon changes discontinuously. We will
consider this possibility further in Section \ref{sec:small-e}. 

However, it's also possible that analyzing the best (2,1)-symmetric tripodal graphon for very small values of $e$ 
is moot. Maybe the best
graphons are tripodal with no symmetry at all, or 4-podal, or 5-podal, and so on. Our data does not rule this out,
nor does it rule out an infinite cascade of phases as we approach $(e,t)=(0,0)$, just as there are infinitely many 
phases with $(n,0)$ or $(n,1)$ or $(n,2)$ symmetry as we approach $(e,t)=(1,1)$. Exploring these 
possibilities is a problem for future research. 

\begin{remark}
Finding the best graphon for $t$ close to $e^3$ is closely connected to the problem of moderate deviations. 
Graphs with $n$ vertices 
and a fixed density $e$ of edges have an expected density of triangles that is close to $e^3$ and a standard deviation
that scales as $n^{-3/2}$. Fluctuations in the triangle density of order $n^{-3/2}$ are called 
{\em small deviations} and are governed by the central limit theorem. Fluctuations of order 1 are called 
{\em large deviations} and are governed by graphons that maximize $S(g)$ subject to constraints on $\edens(g)$
and $\tdens(g)$. Fluctuations of order $n^\alpha$, where $-3/2 < \alpha < 0$, are called {\em moderate deviations}
and are more complicated.

When $e>1/2$, moderate deviations with $\alpha > -1$ and with fewer triangles than expected (``undersaturated graphs'' or ``lower tails'') turn out to be governed by the same graphons that describe the 
$t \to e^3$ limit of large deviations \cite{ADG, NRS2}, while
those with $\alpha < -1$ are qualitatively similar to small deviations \cite{Gold}. However, 
understanding moderate deviations with 
$e<1/2$ is complicated by the fact that we don't know what the optimal graphons are for large deviations. 
That's what we're trying to figure out in this paper! It appears likely that the limiting form of the optimal
graphons as $\delta \to 0$ will shed light on moderate deviations. 

In particular, we can study moderate deviations in the limit as $e \to 0$. When considering upper tails of triangle
counts in sparse graphs, 
the structure derived by \cite{HMS} can be viewed as an $e\to 0$ limit of the graphons found in \cite{KRRS3}. While
we do not know what the best graphon actually is, the ansatz (\ref{eq:ansatz}) with $A=\frac12$ and 
$B=\frac12 - e$ gives a lower bound on the entropy that is much higher than anything previously known. This
suggests a technique for studying moderate deviations for undersaturated sparse graphs. 
\end{remark}

\section{Moving beyond the ansatz} \label{sec:D} We now consider arbitrary tripodal graphons with (2,1)-symmetry. 
To describe such graphons, we need some notation. Let $g_{ij}$ be the value of a graphon on $I_i \times I_j$ and 
let $c_i$ be the width of the $i$-th pode. Because of our (2,1)-symmetry, we have $c_1=c_2$, $g_{11}=g_{22}$ and 
$g_{13}=g_{23}$. The entire graphon can be expressed in terms of four parameters $A$, $B$, $c$, and $D$, with
\begin{eqnarray} \label{eq:21-true}
c_1 = c_2 &=& c/2, \cr c_3 &=& 1-c, \cr 
g_{11} = g_{22} &=& e-A+(1-c)(B+D), \cr 
g_{12} &=& e+A+(1-c)(B+D), \cr 
g_{13}=g_{23} & = & e - cB + \frac{1-2c}{2}D, \cr 
g_{33} & = & e + \frac{c^2}{1-c} B -cD.
\end{eqnarray}
When $D=0$, this reduces to the ansatz (\ref{eq:ansatz}). The extra terms account for changes in the degree
function from one pode to another. The degree function is then 
\be d(x) = \begin{cases} e + \frac{1-c}{2}D & x < c, \cr 
e - \frac{c}{2}D & x>c. \end{cases}
\ee

This graphon can be described more succinctly in terms of the functions
\be v_1(x) = \begin{cases} 1 & x < c/2, \cr -1 & c/2 < x < c, \cr 0 & x>c, \end{cases}
\qquad 
v_2(x) = \begin{cases}\sqrt{1-c} & x<c, \cr \frac{-c}{\sqrt{1-c}} & x>c. \end{cases} \ee
These functions are orthogonal in $L^2([0,1])$, with $\|v_1\|_{L^2}^2=\|v_2\|_{L^2}^2=c$. 
Our graphon is then
\be g(x,y) = e - A v_1(x) v_1(y) + B v_2(x) v_2(y) + \frac{D}{2} \sqrt{1-c} (v_2(x) + v_2(y)). \ee

Since the functions $v_1$ and $v_2$ integrate to zero, our edge density is exactly $e$. Defining $\delta g(x,y) = g(x,y)-e$, the triangle density is easily computed from the expansion 
\begin{eqnarray} 
\tdens(g) &=& e^3 + 3e \int_0^1 (d(x)-e)^2 \, dx + \tdens(\delta g) \cr 
& = & e^3 + \frac34 ec(1-c) D^2 + \frac34 c^2(1-c) BD^2 + c^3(B^3-A^3). 
\end{eqnarray} 

Finally, the entropy is 
\begin{eqnarray}
S & = & \frac{c^2}{2} \Big ( H\big (e-A+(1-c)(B+D)\big) + H\big (e+A+(1-c)(B+D)\big ) \Big ) \cr 
&& + 2c(1-c) H\left (e-cB + \frac{1-2c}{2}D\right ) + (1-c)^2 H\left (e + \frac{c^2}{1-c}B-cD \right ).
\end{eqnarray}
We also compute the first two derivatives of $S$ with respect to $D$ at $D=0$: 
\begin{eqnarray}
\frac{\partial S}{\partial D} & = & \frac12 c^2(1-c) \big ( H'(e-A+B(1-c)) + H'(e+A+B(1-c)\big ) \cr 
&& + c(1-c)(1-c2) H'(e-cB) - c(1-c)^2 H'\left ( e + \frac{c^2}{1-c}B\right ) \cr 
& \approx & \frac{c^2(1-c)}{2} \left ( H'(e\!-\!A\!+\!B(1\!-\!c)) + H'(e\!+\!A\!+\!B(1\!-\!c)) -2H'(e) -2B(1\!-\!c)H''(e) \right ), \cr 
\frac{\partial^2S}{\partial D^2} & = & \frac{c^2(1-c)^2}{2} \left ( H''(e-A+B(1-c)) + H''(e+A+B(1-c)) \right ) \cr 
&& + \frac{c(1-c)(1-2c)^2}{2} H''(e-cB) + c^2(1-c)^2 H''\left ( e+\frac{c^2}{1-c}B \right ) \cr 
& = & \frac{c}{2} H''(e) + O(c^2),
\end{eqnarray} 
where in estimating $\partial S/\partial D$ we have used the linear approximations 
$H'(e-cB) \approx H'(e) -cB H''(e)$ and $H'(e + \frac{c^2}{1-c}B) \approx H'(e) + \frac{c^2}{1-c} B H''(e)$. 

Having $D$ nonzero can increase the entropy by $O(\delta^2D)$, but at a cost of increasing 
the triangle density by $\frac34 c(1-c)(e+Bc) D^2$. 
These must be compensated with other changes to the graphon. Since $\partial S/\partial t$ scales as 
$\delta^{-1}$, and since $c$ scales as $\delta^1$, there is effectively a $\Theta(D^2)$ entropy cost for having
$D \ne 0$, in addition to the $O(cD^2)$ cost from $\partial^2 S/ \partial D^2$. Balancing the cost and benefit, 
the optimal value of $D$ is $O(\delta^2)$, with the result that the $D$ terms only contribute to the entropy at order
$\delta^4$. 

We conclude that the difference between the ansatz graphon (\ref{eq:ansatz}) and the most general
(2,1)-symmetric graphon (\ref{eq:21-true}) is small, especially close the the \ER{} curve. 
It has a non-negligible effect on 
establishing the boundary between the tripodal phase and the symmetric bipodal phase, but does not change the 
qualitative picture. When it comes to a qualitative understanding of the tripodal phase, it is sufficient to work
with the ansatz (\ref{eq:ansatz}). 

\section{The extent of the tripodal phase} \label{sec:limbo}

\begin{figure}[ht]
\includegraphics[width=5in]{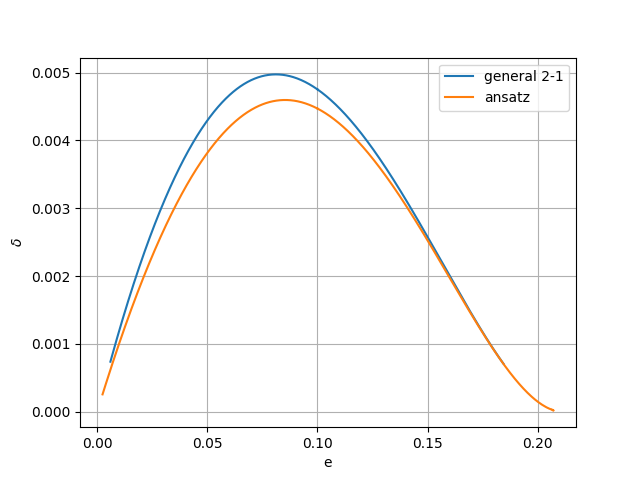}
\caption{The maximum value of $\delta$ for which the ansatz (\ref{eq:ansatz}) or the most general (2,1)-symmetric 
graphon has more entropy than a symmetric bipodal graphon. } \label{fig:limbo}
\end{figure} 

So far, we have only considered graphons with $t$ infinitesimally close to $e^3$. That is, setting $t = e^3-\delta^3$,
we have been looking at the $\delta \to 0$ limit. In this section we consider how big $\delta$ can get and still have
a (2,1)-symmetric tripodal graphon with more entropy than the symmetric bipodal graphon. 
In other words, we want to understand how far down the
(2,1)-symmetric tripodal phase extends in the $(e,t)$ plane. We also examine the nature of the phase transition from
tripodal to symmetric bipodal. 

\subsection{Results} 

The maximum value of $\delta$ is shown in Figure \ref{fig:limbo}, both allowing $D \ne 0$ (the blue curve) and 
restricting our attention to $D=0$ (the orange curve). Note that
\begin{itemize}
\item The largest possible values of $\delta$ in the tripodal phase 
occur at $e \approx 0.08$. Both as $e \to e_0$ and as $e \to 0$, 
the maximum value of $\delta$ goes to zero. 
\item When $e>0.0024$, $\delta$ is never greater than $0.11 e$. Since $t=e^3-\delta^3$, 
this means that tripodal 
phase is almost completely contained in the region $0.998 e^3 < t < e^3$. In terms of the variables $(e,t)$, 
the tripodal phase is extremely small and very easy to miss. 
\item The blue and orange curves are qualitatively similar, but have noticeable differences when $e < 0.12$ or so. 
\item We have calculated the orange curve for essentially the entire interval $[0.0024,e_0]$. However, calculating the 
blue curve is numerically less stable, so we only have results for $e>0.01$. We will discuss the small $e$ regime
further in Section \ref{sec:small-e}.
\end{itemize}

\begin{figure}[ht]
\includegraphics[width=5in]{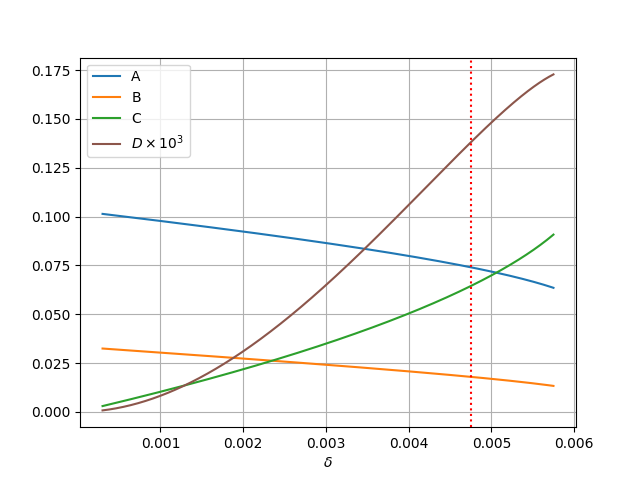} 
\caption{For a fixed value of $e = 0.1$, the optimal values of $A, B, c, D$ as functions of $\delta$, scaled to fit in the same plot. The dotted line shows where the optimal graphon is overtaken in entropy by a symmetric bipodal graphon.}\label{fig:D01}
\end{figure} 
\begin{figure}[ht]
\includegraphics[width=3in]{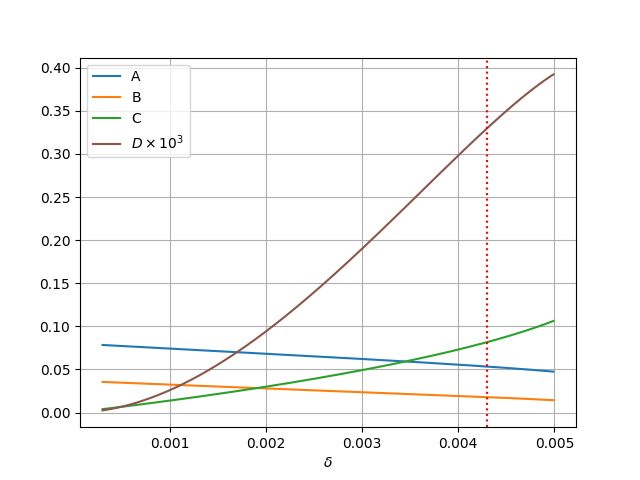} 
\includegraphics[width=3in]{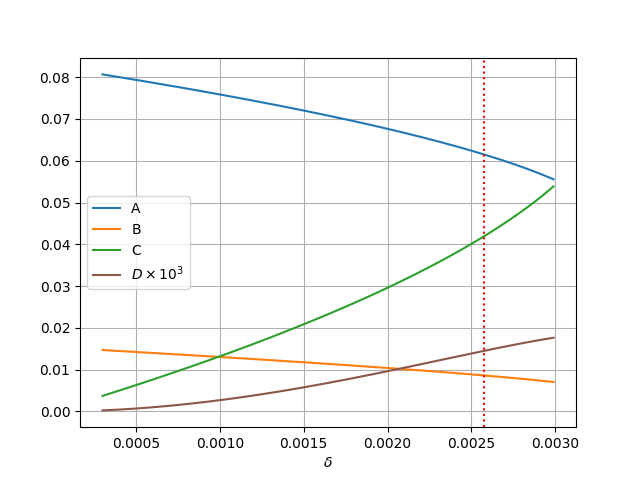} 
\caption{The dependence of $(A, B, c, D)$ on $\delta$ for $e=0.05$ and $e = 0.15$.}\label{fig:D015}
\end{figure} 
\begin{figure}[ht]
\includegraphics[width=5in]{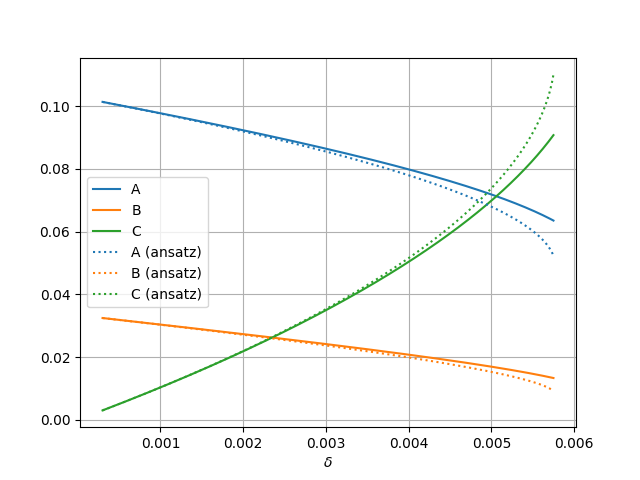} 
\caption{When $e = 0.1$, the optimal values of $(A, B, c)$ for $D=0$ (dotted lines) are nearly identical to those
where we allow $D \ne 0$.}\label{fig:ABc}
\end{figure} 

We next consider how the optimal graphon evolves as we traverse the phase by increasing $\delta$, starting at $\delta =0$. For $e=0.1$, this is shown in Figure \ref{fig:D01}. 
The boundary of the phase at $\delta \approx 0.0048$ is 
not characterized by any dramatic behavior. It is simply the point where the optimal (2,1)-symmetric graphon happens to
have the same entropy as a symmetric bipodal graphon, meaning that the constrained optimization problem does not 
have a unique solution. To the left of the dotted line, the data indicate that the optimal graphon is 
(2,1)-symmetric with the parameter values shown in the figure. To the right of the dotted line, the figure only shows
a local maximum of the entropy, with the actual optimal graphon being symmetric bipodal.

The picture is similar for other values of $e$. See Figure \ref{fig:D015}. 
As $\delta$ increases from zero, $c$ increases linearly and $D$ increases, first
quadratically and then roughly linearly, while $A$ and $B$ gradually decrease. The phase transition does not correspond
to any singularities in the parameters $(A, B, c, D)$. Rather, it is a point where we jump discontinously from a 
graphon described by these parameters to a symmetric bipodal graphon. The
relative size of $D$ versus $A$, $B$ and $c$ is different for different values of $e$, but the graphs of $(A, B, c)$ 
are all qualitatively the same.

Finally, we compare the optimal (2,1)-symmetric graphon to the ansatz (\ref{eq:ansatz}). Figure \ref{fig:ABc} 
shows the situation when $e=0.1$, with the optimal values of $(A,B,c)$ for general $D$ (in solid lines) compared
to the optimal values for $D=0$ (in dotted lines). The differences are minor. Allowing $D \ne 0$ does yield greater 
numerical accuracy, but it doesn't change the overall story. 

\subsection{Methods}

In generating Figure \ref{fig:limbo}, we sought the value of $\delta$ for which the entropy of the best 
possible (2,1)-symmetric graphon equals 
the entropy of a symmetric bipodal graphon. We found this value of $\delta$, which we call $\delta_m$, in three steps. We 
first determined $\delta_m$ to within 0.0001 for $e=0.1$. Next we used Newton's method to home in on the true value of 
$\delta_m$ for $e=0.1$. Finally, we varied $e$ and used Newton's method to track the changes in $\delta_m$ as a function
of $e$. We did this 3-step process twice, once assuming $D=0$ (that is, the ansatz  (\ref{eq:ansatz})) and once allowing
$D$ to float. 

We now describe the three steps for the ansatz and then note the differences for $D \ne 0$. 
\begin{enumerate}

\item We began at $(e, \delta) = (0.1, 0)$ and chose $A$ and $B$ to maximize $F(A,B)$, 
as in Section \ref{sec:F}. 
We then incremented $\delta$ in steps of $0.0001$. At each new value of $\delta$, we solved for
$c$ in terms of $(A,B)$ and wrote the entropy as a function of $(A,B)$. We maximized that entropy, using the optimal 
values of $(A,B)$ from the previous value of $\delta$ as the starting point for Newton's method. We recorded the resulting
entropy and compared it to that of a symmetric bipodal graphon, continuing until the tripodal entropy dropped below the
symmetric bipodal entropy. The true value of $\delta_m$ then lay between the last two values of $\delta$.  

\item Once we knew $\delta_m$ to within 0.0001 and knew the corresponding values of $(A,B)$, we alternated between two forms of
Newton's method to obtain more accuracy. First we held $(A,B)$ fixed and varied $\delta$ to equate the entropy of our 
tripodal graphon with the entropy of a symmetric bipodal graphon. Then we held $\delta$ fixed and varied $(A,B)$ to maximize
the entropy for the given $\delta$. We alternated between these two calculations, converging to a triple $(A,B,\delta)$
that maximized entropy for the given $\delta$ and that had the same entropy as a symmetric bipodal graphon. This gave us
$\delta_m$ for $e=0.1$ to machine accuracy. 

\item For other values of $e$, we did not start at $\delta=0$. Instead, we incremented (or decremented) $e$ in steps of 0.001, 
using the values
of $(A, B, \delta_m)$ from the previous value of $e$ as the starting point for an alternating 
pair of Newton's method calculations, as in Step 2. 
This filled out the rest of the orange curve in Figure \ref{fig:limbo}. 

\end{enumerate}

The calculations for $D \ne 0$ were similar, only with a few adjustments. Since the optimal $D$ is $O(\delta^2)$, 
we started at $(e,\delta)=(0.1,0)$, 
with the same values of $(A,B)$ as before and with $D=0$. For each new $\delta$, we eliminated $A$ rather 
than $c$ and used Newton's 
method to obtain values of $(B, c, D)$ that maximized the entropy. In the second and third steps, we alternated between
varying $\delta$ while holding $(B,c,D)$ fixed (to equate the tripodal and bipodal entropies) and varying $(B,c,D)$ with 
$\delta$ fixed (to maximize the tripodal entropy). These calculations generated the blue curve. 

To generate Figures \ref{fig:D01}--\ref{fig:ABc}, we repeated the first step of the above calculation. That is, we started at 
$(e,\delta) = (0.1,0)$ or $(0.05,0)$ or $(0.15,0)$, with $D=0$ and with $(A,B)$ maximizing $F(A,B)$. 
We gradually increased
$\delta$, eliminated one variable ($c$ or $A$), and used several rounds of Newton's method on the remaining variables 
to maximize the entropy,
using data from the previous value of $\delta$ as a starting point. 
This gave $(A,B,c)$ or $(A,B,c,D)$ as functions of $\delta$. 

\section{Graphs with low edge density}\label{sec:small-e}

We now examine graphons whose edge density is close to zero. As we saw in Section \ref{sec:F}, there are two
local maximizers of $F(A,B)$, one with $A \approx 2.5e$ and $B \approx 1.5e$ and the other with $A \approx \frac12$
and $B \approx \frac12 - e$. As $\delta$ increases, this gives rise to two families of (2,1)-symmetric tripodal graphons,
each of which is a local entropy maximizer in the space of (2,1)-symmetric graphons. Based on the size of $A$ and $B$,
we refer to the first as the ``$O(e)$ solution'' and the second as the ``$\Theta(1)$ solution''. 

In Section \ref{sec:limbo} we compared the entropy of the $O(e)$ solution to the entropy of a symmetric bipodal graphon
and determined the largest value of $\delta$ for which the first was greater than the second. Here we compare the 
entropies of the $O(e)$ and $\Theta(1)$ solutions to each other and to that of a symmetric bipodal graphon, using essentially the 
same numerical methods as in Section \ref{sec:limbo}. 

\begin{figure}[ht]
\includegraphics[width=4in]{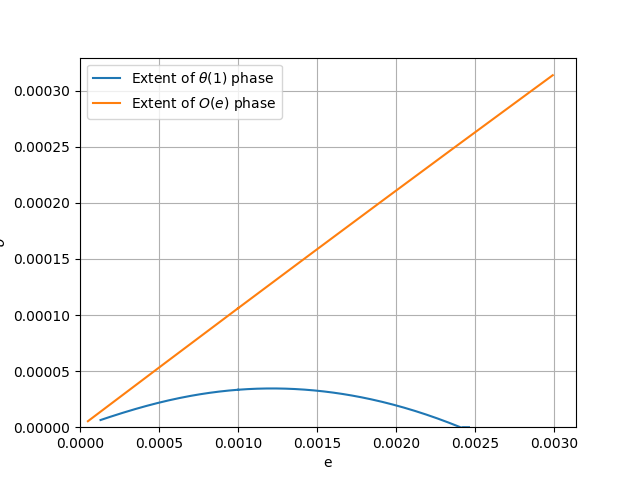} 
\caption{As a function of $e$, the values of $\delta$ at which the $\Theta(1)$ phase transitions to $O(e)$ (blue) and at which the $O(e)$ phase transitions to symmetric bipodal (orange).}\label{fig:small-e}
\end{figure} 

The results are summarized as follows:
\begin{itemize}
\item When $e>0.0024$ or so, the $\Theta(1)$ solution never does better than the $O(e)$ solution. 
The entire tripodal
region takes the $O(e)$ form and extends from \ER{} to the value of $\delta$ given by the orange curve
of Figure \ref{fig:small-e} or the blue curve of Figure \ref{fig:limbo}. This is the
phase that we previously studied in Section \ref{sec:limbo}.
\item When $e < 0.0024$, $F(A,B)$ is maximized with $(A,B)=\Theta(1)$, so the optimal graphon for sufficiently
small values of $\delta$ takes the $\Theta(1)$ form. However, the $O(e)$ solution does better for
larger values of $\delta$, outperforming both $\Theta(1)$ and symmetric bipodal. 
See Figure \ref{fig:small-e} for the values of $\delta$ where we transition from $\Theta(1)$ to $O(e)$, and 
where we then transition from $O(e)$ to symmetric bipodal.
\item Considering all values of $e$, the tripodal region consists of two distinct phases. The 
$\Theta(1)$ phase hugs the \ER{} curve very tightly for $e<0.0024$ and is 
separated from the symmetric bipodal phase
by the $O(e)$ phase. This phase is similar in shape to the $O(e)$ phase, which hugs the ER curve for 
$e < e_0 \approx 0.2113$. However, the $\Theta(1)$ phase is {\em much} smaller, being limited to 
$e<0.0024$ and roughly $0.9999 e^3 < t < e^3$. 
\item When speaking of the $O(e)$ and $\Theta(1)$ phases, we have implicitly assumed that the optimal 
graphon is either bipodal or 
(2,1)-symmetric tripodal. It is possible that part of what we call the $O(e)$ phase, and part or all of the $\Theta(1)$
phase, is actually in an asymmetric tripodal phase or a more complicated multipodal phase. We have seen no
evidence to suggest such complicated structures, but at this point we cannot rule them out. 
\end{itemize}

\subsection*{Acknowledgments} We thank Charles Radin for helpful discussions. This work is partially supported
by the National Science Foundation. For access to our numerical data, please contact the first author.

\end{document}